\documentclass{article}
\input{tcilatex}
\begin{document}

\begin{center}
\bigskip \textbf{Some Characterizations for Curves by the Help Of Spherical
Representations in the Galilean and Pseudo-Galilean Space}$.$

\textbf{\bigskip }

{\small Alper Osman \"{O}\u{G}RENM\.{I}\c{S}, M\"{u}nevver YILDIRIM YILMAZ,
Mihriban K\"{U}LAHCI}

{\small F\i rat University, Science Faculty, Mathematics Department}

{\small 23119 Elaz\i \u{g} / T\"{U}RK\.{I}YE}

{\small ogrenmisalper@gmail.com, myildirim@firat.edu.tr,
mihribankulahci@gmail.com}

\bigskip
\end{center}

\textbf{Abstract. }In this paper, we focus on some characterizations for
curves in the Galilean and Pseudo-Galilean\textbf{\ }space.

\textbf{Keywords. }Galilean space, pseudo-Galilean space, frenet frame,
spherical representation, harmonic curvature, .

\textbf{MSC(2000). }53A35; 53B30.

\bigskip

\textbf{1. Introduction. }

Discovering Galilean space-time is probably one of the major achievements of
nonrelativistic physics. One may consider Galilean space as the limit case
of a pseudo-Euclidean space in which the isotropic cone degenerates to a
plane. This limit transition corresponds to the limit transition from the
special theory of relativity to classical mechanics. Nowadays Galilean space
is becoming increasingly popular as evidenced from the connection of the
fundamental concepts such as velocity, momentum, kinetic energy, etc. and
principles as indicated in [8].

As a branch of science, geometry of space is associated with mathematical
group. The idea of invariance of geometry under transformation group shows
that, on some spacetimes of maximum symmetries there should be a principle
of relativity, which requires the invariance of physical laws without
gravity under transformations among inertial systems.

Galilean space G$_{3}$ has been investigated in [3], [4], [5], [6], [7] from
the differential geometrical point of view.The mathematical model of
Galilean space is a three dimensional complex projective space P$_{3}$ in
which the absolute figure \{ $w$, $f$, $I_{1},$ $I_{2}$ \} consists of a
real plane $w$ (the absolute plane), a real line $f\subset w$ (the absolute
line) and two complex conjugate points $I_{1},$ $I_{2}\in f$ (the absolute
points).

We may also take, as a real model of the space $G_{3},$ a real projective
space $P_{3}$ with the absolute \{ $w,$ $f$ \} consisting of a real plane $%
w\subset G_{3}$ and a real line $f\subset w$ on which an elliptic involution 
$\varepsilon $ has been defined.

Using homogeneous coordinates one may write

\begin{eqnarray*}
w...x_{0} &=&0,~\text{~~}f...x_{0}=x_{1}=0 \\
\varepsilon &:&(0:0:x_{2}:x_{3})\rightarrow (0:0:x_{3}:-x_{2})
\end{eqnarray*}

\bigskip On the contrary, for nonhomogeneous coordinates of the similarity
group $H_{8}$ has the following form,

\begin{eqnarray}
x^{\prime } &=&a_{11}+a_{12}x,  \nonumber \\
y^{\prime } &=&a_{21}+a_{22}x+a_{23}\cos \varphi y+a_{23}\sin \varphi z, 
\TCItag{1} \\
z^{\prime } &=&a_{31}+a_{32}x-a_{23}\sin \varphi y+a_{23}\cos \varphi z, 
\nonumber
\end{eqnarray}%
here we denote $a_{ij}$ and $\varphi $ as real numbers.

Taking $a_{12}=a_{23}=1$ we obtain the subgroup $B_{6}-$ the group of
Galilean motions as follows;

\begin{eqnarray*}
x^{\prime } &=&a+x \\
B_{6}...y^{\prime } &=&b+cx+y\cos \varphi +z\sin \varphi \\
z^{\prime } &=&d+ex-y\sin \varphi +z\cos \varphi .
\end{eqnarray*}

We may observe four classes for lines in $G_{3}$ as indicated below:

a) (proper) nonisotropic lines-they don't meet the absolute line $f.$

b) (proper) isotropic lines -lines that don't belong to the plane $w$ but
meet the absolute line $f.$

c) unproper nonisotropic lines-all lines of $w$ but $f.$

d) the absolute line $f.$

Planes $x=const.$ are Euclidean and so is the plane $w.$ Other planes are
isotropic.

Here the coefficients $a_{12}$ and $a_{23}$ play the special role.

In particular, for $a_{12}=a_{23}=1$ (1) shows the group $B_{6}\subset H_{8}$
of isometries of the Galilean space $G_{3}$ [6].

\bigskip

\textbf{2. Basic notions and properties}

Let $\alpha :I\rightarrow G_{3},$ $I\subset IR$ be a curve given by

\[
\alpha (t)=(x(t),y(t),z(t)), 
\]%
where $x(t),$ $y(t),$ $z(t)$ $\in C^{3}$ (the set of three times
continuously differentiable functions) and $t$ run through a real interval
[6].

Let $\alpha $ be a curve in $G_{3},$ parameterized by arc length $t=s,$
given in coordinate form

\begin{equation}
\alpha (s)=(s,y(s),z(s)).  \tag{2}
\end{equation}

Then the curvature $\kappa (s)$ and the torsion $\tau (s)$ are defined by

\begin{eqnarray}
\kappa (s) &=&\sqrt{y^{^{\prime \prime }2}(s)+z^{^{\prime \prime }2}(s)} 
\TCItag{3} \\
\tau (s) &=&\frac{\det (\alpha ^{^{\prime }}(s),\alpha ^{^{^{\prime \prime
}}}(s),\alpha ^{^{^{\prime \prime \prime }}}(s))}{\kappa ^{2}}  \nonumber
\end{eqnarray}%
and associated moving trihedron is given by

\begin{eqnarray}
T(s) &=&\alpha ^{\prime }(s)=(1,y^{^{\prime }}(s),z^{^{\prime }}(s)) 
\TCItag{4} \\
N(s) &=&\frac{1}{\kappa (s)}\alpha ^{^{\prime \prime }}(s)=\frac{1}{\kappa
(s)}(0,y^{^{^{\prime \prime }}}(s),z^{^{^{\prime \prime }}}(s))  \nonumber \\
B(s) &=&\frac{1}{\kappa (s)}(0,-z^{^{^{\prime \prime }}}(s),y^{^{^{\prime
\prime }}}(s)).  \nonumber
\end{eqnarray}

The vectors $T,$ $N,$ $B$ are called the vectors of the tangent, principal
normal and binormal line of $\alpha ,$ respectively. For their derivatives
the following Frenet formulas hold

\begin{eqnarray}
T^{^{\prime }} &=&\kappa N  \TCItag{5} \\
N^{^{\prime }} &=&\tau B  \nonumber \\
B^{^{\prime }} &=&-\tau N  \nonumber
\end{eqnarray}

Scalar product in the Galilean space $G_{3}$ is defined by

\begin{equation}
g(X,Y)=\{%
\begin{array}{cc}
x_{1}y_{1}, & if\text{ \ }x_{1}\neq 0\text{ }\vee y_{1}\neq 0 \\ 
x_{2}y_{2}+x_{3}y_{3}, & if\text{ \ }x_{1}=0\text{ }\wedge \text{ }y_{1}=0%
\end{array}%
\text{ ,}  \tag{6}
\end{equation}%
where \ $X=(x_{1},x_{2},x_{3})$ and $Y=(y_{1},y_{2},y_{3})$

\textbf{Definition 2.1. }Let $\alpha $ be a curve in 3-dimensional Galilean
space $G_{3}$, and $\{T,N,B\}$ be the Frenet frame in 3-dimensional Galilean
space $G_{3}$ along $\alpha $ . If $\kappa $ and $\tau $ are positive
constants along $\alpha $ , then $\alpha $ is called a circular helix with
respect to the Frenet frame [4]

\textbf{Definition 2.2.} .Let $\alpha $ be a curve in 3-dimensional Galilean
space $G_{3}$, and $\{T,N,B\}$ be the Frenet frame in 3-dimensional Galilean
space $G_{3}$ along $\alpha $ . A curve $\alpha $ such that%
\[
\frac{\kappa }{\tau }=const. 
\]%
is called a general helix with respect to Frenet frame [4].

\textbf{Definition 2.3. }Let's give curve of $\alpha \subset G_{3}$ with
coordinate neighbourhood $(I,\alpha ).$ $\kappa (s)$ and $\tau (s)$ be a
curvature of $\alpha $ on point of $\alpha (s)\in \alpha $ corresponding to $%
\forall s\in I.$ The function of $H$ can be defined as

\begin{eqnarray*}
H &:&I\rightarrow R \\
s &\rightarrow &H(s)=\frac{\kappa }{\tau }
\end{eqnarray*}%
where the function of $H$ is called $1^{st}$ harmonic curvature of $\alpha $
on point of $\alpha (s).$

\textbf{Remark 2.1. } Similar definitions can be given in the
pseudo-Galilean space.

\textbf{Remark 2.2.} In [2], for the pseudo-Galilean Frenet trihedron of an
\linebreak admissible curve $\alpha ,$ the following derivative Frenet
formulas are true.

\begin{eqnarray*}
T^{^{\prime }}(s) &=&\kappa (s)N(s) \\
N^{^{\prime }}(s) &=&\tau (s)B(s) \\
B^{^{\prime }}(s) &=&\tau (s)N(s)
\end{eqnarray*}%
where $T(s)$ is a spacelike, $N(s)$ is a spacelike and $B(s)$ is a timelike
vector, $\kappa (s)$ is the pseudo-Galilean curvature given by above
equations and $\tau (s)$ is the pseudo-Galilean torsion of $\alpha $ defined
by

\[
\tau (s)=\frac{y^{^{\prime \prime }}(s)z^{^{\prime \prime \prime
}}(s)-y^{^{\prime \prime \prime }}(s)z^{^{\prime \prime }}(s)}{\kappa ^{2}(s)%
}. 
\]

\bigskip

\textbf{3. The Arc Length Of Spherical Representations Of The Curve }$\alpha
\subset G_{3}$

In this section, using method in [1], some characterizations related to
spherical representations are obtained in Galilean and Pseudo-Galilean
3-space.

\textbf{Theorem 3.1. }$\alpha \subset G_{3}$ is an ordinary helix if and
only if

\[
s_{_{T}}=\tau Hs+c. 
\]

\textbf{Proof. }Let \ $T=T(s)$ be the tangent vector field of the curve

\begin{eqnarray*}
\alpha &:&I\subset R\rightarrow G_{3} \\
s &\rightarrow &\alpha (s)
\end{eqnarray*}

The spherical curve $\alpha _{T}=T$ on $S^{2}$ is called first spherical
representation of the tangent of $\alpha .$

Let $s$ be the arc length parameter of $\alpha .$ If we denote the arc
length of the curve $\alpha _{_{T\text{ \ }}}$by $s_{_{T}}$, then we may
write 
\[
\alpha _{_{T}}(s_{_{T}})=T(s). 
\]

Letting $\frac{d\alpha _{_{T}}}{ds_{_{T}}}=T_{_{T}}$ we have $%
T_{_{T}}=\kappa \overrightarrow{N}\frac{ds}{ds_{_{T}}}.$ Hence we obtain $%
\frac{ds_{_{T}}}{ds}=\kappa .$ Thus we give the following result.

If $\kappa $ is the first curvature of the curve $\alpha :I\rightarrow
G_{3}, $ then the arc length $s_{_{T}}$ of the tangentian representation $%
\alpha _{_{T}}$ of $\alpha $ is 
\[
s_{_{T}}=\dint \kappa ds+c. 
\]

If the harmonic curvature of $\alpha $ is $H=\frac{\kappa }{\tau },$ we get

\[
ds_{_{T}}=\dint \tau Hds+c 
\]%
where \ $c$ is an integral constant.

\textbf{Theorem 3.2. }$\alpha \subset G_{3}$ is an ordinary helix if and
only if

\[
s_{_{N}}=\frac{\kappa }{H}s+c. 
\]

\textbf{Proof. }Let \ $\overrightarrow{N}=\overrightarrow{N}(s)$ be the
principal normal vector field of the curve

\begin{eqnarray*}
\alpha &:&I\subset R\rightarrow G_{3} \\
s &\rightarrow &\alpha (s)
\end{eqnarray*}

The spherical curve $\alpha _{_{N}}=\overrightarrow{N}$ on $S^{2}$ is called
second spherical representation for $\alpha $ or is called the spherical
representation of the principal normals of $\alpha .$ Let $s\in I$ be the
arc length parameter of $\alpha .$ If we denote the arc length of the curve $%
\alpha _{_{N\text{ \ }}}$by $s_{_{N}}$, we may write 
\[
\alpha _{_{N}}(s_{_{N}})=\overrightarrow{N}(s). 
\]

Moreover letting $\frac{d\alpha _{_{N}}}{ds_{_{N}}}=T_{_{N}}$ we obtain

\[
T_{N}=\tau \overrightarrow{B}\frac{ds}{ds_{_{N}}}. 
\]

Hence we have

\[
\frac{ds_{_{N}}}{ds}=\tau . 
\]

Thus we give the following result.

If $\tau $ is the second curvature of the curve $\alpha :I\rightarrow G_{3},$
then the arc length $s_{_{N}}$ of the principal normal representation $%
\alpha _{_{N}}$ of $\alpha $ is 
\[
s_{_{N}}=\dint \tau ds+c. 
\]

If the harmonic curvature of $\alpha $ is $H=\frac{\kappa }{\tau },$ we get

\[
s_{_{N}}=\dint \tau ds+c. 
\]

If the harmonic curvature of $\alpha $ is $H=\frac{\kappa }{\tau },$ we get

\[
s_{_{N}}=\dint \frac{\kappa }{H}ds+c 
\]

where \ $c$ is an integral constant.

\textbf{Theorem 3.3. }$\alpha \subset G_{3}$ is an ordinary helix if and
only if

\[
s_{_{B}}=\frac{\kappa }{H}s+c. 
\]

\textbf{Proof. }Let \ $\overrightarrow{B}=\overrightarrow{B}(s)$ be the
binormal vector field of the curve

\begin{eqnarray*}
\alpha &:&I\subset R\rightarrow G_{3} \\
s &\rightarrow &\alpha (s)
\end{eqnarray*}

The spherical curve $\alpha _{_{B}}=\overrightarrow{B}$ on $S^{2}$ is called
third spherical representation for $\alpha $ and the spherical
representation of the binormal of $\alpha .$

Let $s\in I$ be the arc length parameter of $\alpha .$ If we denote the arc
length parameter of the curve $\alpha _{_{B\text{ \ }}}$by $s_{_{B}}$, we
may write 
\[
\alpha _{_{B}}(s_{_{B}})=\overrightarrow{B}(s). 
\]

Moreover letting $\frac{d\alpha _{_{B}}}{ds_{_{B}}}=T_{_{B}},$ we obtain

\[
T_{B}=-\tau \overrightarrow{N}\frac{ds}{ds_{B}}. 
\]

Hence we have $\frac{ds_{_{B}}}{ds}=\tau $ and $s_{_{B}}=\dint \tau ds+c$ or
in terms of the harmonic curvature of $\alpha $ we obtain

\[
s_{_{B}}=\dint \frac{\kappa }{H}ds+c. 
\]

\textbf{Note. }Same theorems can be given in the Pseudo-Galilean 3-space.

\bigskip

\textbf{References}

[1] Hac\i saliho\u{g}lu, H. H., A New Characterization For Inclined Curves
By The Help Of Spherical Representations, International Electronic Journal
Of Geometry, Vol. 2, No.2 , (2009) 71-75.

[2] Divjak, B., Curves in Pseudo-Galilean Geometry, Annales Univ. Sci.
Budapest, 41: (1998) 117-128.

[3] Kamenarovic, I., Existence Theorems for Ruled Surfaces In The Galilean
Space $G_{3}$ , Rad HAZU Math. 456 (10), (1991) 183-196.

[4] \"{O}\u{g}renmi\c{s}, A. O., Bekta\c{s}, M. and Erg\"{u}t ,M., On The
helices in the Galilean space $G_{3}$, Iranian Journal of science \&
Technology, Transaction A , Vol.31, No:A2, (2007) 177-181.

[5] \"{O}\u{g}renmi\c{s}, A. O., \"{O}ztekin, H., Erg\"{u}t, M., Bertrand
Curves In Galilean Space and Their Characterizations, Kragujevac J.Math. 32
(2009) 139-147.

[6] Pavkovic, B. J., Kamenarovic, I., The Equiform Differential Geometry of
Curves in the Galilean Space, Glasnik Matematicki, 22(42) (1987) 449-457.

[7] R\"{o}schel, O., Die Geometrie Des Galileischen Raumes, Berichte der
Math.-Stat. Sektion im Forschungszentrum Graz, Ber. 256, (1986) 1-20.

[8] Yaglom, I. M., A Simple Non-Euclidean Geometry and Its Physical Basis,
Springer-Verlag, New York Inc. 306, (1979) 201-214.

\end{document}